\documentclass[preprintnumbers,amsmath,amssymb]{revtex4}

\usepackage{graphicx}
\usepackage{bm}

\begin{document}

\title{Periodic  travelling wave solutions of discrete nonlinear Schr\"odinger equations}

\author{Dirk Hennig}
\affiliation{Eberswalde, Germany}

\date{\today}

\begin{abstract}
\noindent The existence of nonzero periodic  travelling wave solutions for 
a general discrete nonlinear Schr\"odinger
equation (DNLS) on finite  one-dimensional lattices is proved. The DNLS features a 
general nonlinear term and variable range of interactions 
going beyond the usual nearest-neighbour interaction. 
The problem of  the existence of travelling wave solutions is 
converted into a fixed point problem  for an operator on some appropriate function space  
which is solved by means of Schauder's Fixed Point Theorem. 
\end{abstract}

\maketitle

\noindent 

Travelling waves (TWs) in lattice dynamical systems have attracted considerable interest recently 
(see e.g. \cite{Zinner}-\cite{Guo}). A variety of systems has been studied including 
Fermi-Pasta-Ulam systems, discrete nonlinear Klein–Gordon systems,  phase oscillator lattices, 
Josephson junction systems, reaction-diffusion systems  and the discrete nonlinear Schr\"odinger equation.
Some exact results 
regarding the existence, stability and uniqueness of TWs in the above mentioned
systems have been obtained.

Here we are interested in the existence   of periodic TW solutions of the 
following general discrete nonlinear Schr\"odinger equation on 
finite one-dimensional lattices
\begin{equation}
i\frac{d \psi_n}{dt}=\sum_{j =1}^{N_c}\,\kappa_j\,[\psi_{n+j}-2\psi_n+\psi_{n-j}]+ 
F(|\psi_n|^2)\psi_n,\,\,\,1 \le n \le N,\label{eq:system}
\end{equation}
with $\psi_n \in {\mathbb{C}}$.

The solutions satisfy periodicity conditions:
\begin{equation}
\psi_{N+m}(t)=\psi_m(t),\label{eq:pcs}
\end{equation}
 for  $m \in {\mathbb{Z}}$, viz. we consider the DNLS on rings.  Each unit interacts with its $N_c$ neighbouring oscillators to 
the left and right respectively.  $N_c=1,...,[(N-1)/2]$ determines the interaction radius, 
which ranges from nearest-neighbour interaction obtained for $N_c=1$ to global coupling 
when $N_c=(N-1)/2$ and $N$ is odd.

\vspace*{0.5cm}

Assume the following condition on $F(|\psi_n|^2)$ holds:

\vspace*{0.5cm}

\noindent {\bf{A:}} 
$F\in C({\mathbb{R}}_+,{\mathbb{R}})$ for ${\mathbb{R}}_+=[0,\infty)$, $F(0)=0$. 
There are constants $a>0$, $b>0$ such that 
\begin{equation}
 |F(x)|<a(1+x^{b}),\label{eq:A1}
\end{equation}
for any $x\ge 0$.

\vspace*{0.5cm}

The standard DNLS, arising for $F(|\psi_n|^2)=|\psi_n|^2$ and $\kappa_1\ne 0$, 
$\kappa_{j \ne 1}= 0$ in (\ref{eq:system}),  is known to support 
periodic travelling wave solutions (see e.g. \cite{TWS}). 
In \cite{DNLS} existence and bifurcation results were derived using variational methods
for periodic travelling waves of the same general DNLS system
as given in (\ref{eq:system})  on the lattice ${\mathbb{Z}}$.
In this note  we present a proof of 
the existence of periodic travelling wave 
solutions of (\ref{eq:system}) on finite lattices.
To obtain our result, some appropriate function space is introduced on which the original 
problem is presented as a fixed point problem for a corresponding operator. 
By exploiting Schauder's 
Fixed Point Theorem the existence of periodic travelling wave solutions is 
established.

The system (\ref{eq:system}) possesses two conserved quantities, the energy
\begin{equation}
 {\cal H}=\sum_{n=1}^N\,\left[G(|\psi_n|^2)-\sum_{m\neq n}\kappa_m|\psi_m-\psi_n|^2\right],\,\,\,G(x)=\int_0^xf(x)dx,
\end{equation}
and the power
\begin{equation}
 {\cal P}=\sum_{n=1}^N\,|\psi_n|^2.\label{eq:norm}
\end{equation}

We consider travelling wave solutions of the form:
\begin{equation}
 \psi_n(t)=\Psi\left(kn-\omega \,t\right),
\end{equation}
with a $2\pi-$periodic function $\Psi(u)$, $u=kn-\omega\,t$, where $k\in(-\pi,\pi)$ and 
$\omega \in {\mathbb{R}}\setminus \{0\}$ are the wave parameters. 
 
In order that a travelling wave solution satisfies the periodicity 
conditions in (\ref{eq:pcs}) we adopt the lattice size accordingly.
This means that for a given wavenumber $|k| =\pi q$ with rational $q=r/s$ 
and two relatively prime integers  $r,s \in {\mathbb{Z}}_+\setminus \{0\}$,  $r<s$, 
the   number of sites of the lattice, $N$, is  supposed to be 
an appropriate multiple of the  the minimal spatial period 
of the associated periodic travelling wave, determined by $L=s/r$, so that the periodicity  
  conditions in (\ref{eq:pcs}) are fulfilled.
  
\vspace*{0.5cm}
 
Regarding the existence of periodic travelling wave solutions we state the following:

\vspace*{0.5cm}
\noindent{\bf Theorem:}{\it \,\,Let {\bf{A:}}  hold.  Then for any rational number  
$q \in {\mathbb{Q}} \cap (0,1)$ 
 there exists nonzero  periodic travelling wave solution 
$\psi_n(t)=\Psi\left(kn-\omega\,t\right)\equiv\Psi(u)$ of (\ref{eq:system}) with 
$\Psi \in C^1(\mathbb{R},\mathbb{C})$, 
such that
\begin{equation}
 \Psi(u+2\pi)=\Psi(u),\,\,\forall u \in {\mathbb{R}},
\end{equation}
provided that
\begin{equation}
|\omega| \ge {\cal R}\left(
1+p\frac{a(1+{\cal P}^b)}{[{\cal R}(1+q)+4{\bar{\kappa}}+
a(1+{\cal {P}}^b)}\right),\label{eq:As}
\end{equation}
where 
\begin{equation}
 \bar{\kappa}=\sum_{j=1}^{N_c}\kappa_j,
\end{equation}
\begin{equation}
 p={\tilde{q}}+\frac{1}{1-q},\label{eq:p}
\end{equation}
\begin{equation}
 {\tilde{q}} =\left\{ \begin{array}{ccl} \frac{1}{q} & \mbox{for}
& 0<q<\frac{1}{2} \\
 \\\frac{1}{1-q} & \mbox{for} & \frac{1}{2}\le q <1,
\end{array}\right.\label{eq:tildeq}
\end{equation}
and 
${\cal{R}}$ determines the range $[- {\cal R}, {\cal R}]$ of the 
function $g\in C({\mathbb{R}}\setminus\{0\},{\mathbb{R}})$ given by
\begin{equation}
g(x)=\frac{2}{x}\sum_{j=1}^{N_c}\kappa_j\,\sin^2(jx).
\end{equation}
}

\vspace*{1.0cm}

In the following we reformulate the original problem as a fixed point
problem in a Banach space in a similar vein to the approach in \cite{SIAM}.

To prove the assertions of the Theorem we utilise Schauder's Fixed Point Theorem (see e.g. in \cite{Zeidler}): 
{\it Let $M$ be a closed convex subset of a 
Banach space $X$. Suppose
$T:\,\,M\rightarrow\,M$  is continuous mapping such that $T(M)$ is a relatively compact subset of $M$.
Then  $T$ has a fixed point.}

\vspace*{1.0cm}

\noindent{\bf Proof:}
 Travelling wave solutions $\Psi$ satisfy the advance-delay equation
\begin{equation}
 -i \omega \Psi^{\prime}(u)=\sum_{j=1}^{N_c} \kappa_j 
 \Delta_j\Psi(u)+F(|\Psi(u)|^2)\Psi(u),\label{eq:TW}
\end{equation}
where  
$\Delta_j\Psi(u)=\Psi(u+j)-2\Psi(u)+\Psi(u-j)$ and 
$\Psi(u+2\pi)=\Psi(u),\,\,\forall u \in {\mathbb{R}}$, so that 
according to the Bloch-Floquet Theorem a solution  must be of the form 
\begin{equation}
 \Psi(u)=\exp(i \,q u)\Phi(u),\label{eq:BFT}
\end{equation}
where $q \in {\mathbb{Q}} \cap (0,1)$ and 
\begin{equation}
 \Phi(u+2\pi)=\Phi(u),\,\,\,\forall u \in {\mathbb{R}}.\label{eq:2pi}
\end{equation}
 
Substituting (\ref{eq:BFT}) into 
the equation (\ref{eq:TW}) one obtains:
\begin{equation}
 -i\omega \Phi^{\prime}(u) +\omega q\,\Phi(u)=\sum_{j=1}^{N_c} \kappa_j 
 {\tilde{\Delta}}_j\Phi(u)+ F(|\Phi(u)|^2)\Phi(u),\label{eq:ref}
\end{equation}
with  
\begin{equation}
 {\tilde{\Delta}}_j\Phi(u)=\Phi(u+j)\exp(iqj)-2\Phi(u)+\Phi(u-j)\exp(-iqj).
\end{equation}

Thus, the task amounts to find   
 $2\pi-$periodic functions $\Phi\in C^{1}({\mathbb{R}},{\mathbb{C}})$) satisfying 
 Eq.\,(\ref{eq:ref}).

\vspace*{1.0cm}

For the forthcoming discussion  Eq.\,(\ref{eq:ref}) is 
suitably re-arranged as follows:
\begin{equation}
 -i\omega \Phi^{\prime}(u) +\omega q\,\Phi(u)-\sum_{j=1}^{N_c} \kappa_j 
 {\tilde{\Delta}}_j\Phi(u)= F(|\Phi(u)|^2)\Phi(u).\label{eq:ref1}
\end{equation}
Note that terms nonlinear in $\Phi$ feature only on the r.h.s. of (\ref{eq:ref1}).

Let $q \in {\mathbb{Q}} \cap (0,1)$ be fixed. We identify ${\mathbb{C}}$ with ${\mathbb{R}}^2$. 
Denote by ${{X}}_q^{h}$ the real Banach spaces
\begin{equation}
 {{X}}_q^h\,=\,\left\{\,\Theta\in C^{\,h}_{2\pi}({\mathbb{R},\mathbb{C}})\,\right\}\,,\,\,h=0,1,
\end{equation}
where $C^{\,h}_{\pi}(\mathbb{R},{\mathbb{C}})$ 
is the Banach space of $2\pi-$periodic and 
$C^{\,h}$ functions $\Theta\,:\,{\mathbb{R}}\,\rightarrow\,{\mathbb{C}}$ 
equipped with 
norms  given by
\begin{equation}
 \parallel \Theta \parallel_{C^{0}_{2\pi}}=
 \max_{u\in [0,2\pi]}|\Theta(u)|\,,\,\,\,\Theta 
 \in C^0_{2\pi}({\mathbb{R},\mathbb{C}}).
\end{equation}
and 
\begin{equation}
 \parallel \Theta \parallel_{C^{1}_{2\pi}}=
 \max_{u\in [0,2\pi]}|\Theta(u)|+\max_{u\in (0,2\pi)}|\Theta^{\prime}(u)|\,,\,\,\,\Theta 
 \in C^1_{2\pi}({\mathbb{R},\mathbb{C}}),
\end{equation} 
respectively. $X_q^1$ is compactly embedded in $X_q^0$ ($X_q^1\Subset X_q^{0}$).

We decompose functions $\Theta\in X_q^{1}$  in a Fourier series  
\begin{equation}
 \Theta(u)=\sum_{l\in {\mathbb{Z}}}\Theta_l \exp(i\,lu).\label{eq:Fourier}
\end{equation}

Related with the l.h.s. of (\ref{eq:ref1}) we consider the linear mapping: $M_q\,:\,X^1_q\,\rightarrow\,X^{0}_q$:
\begin{equation}
 M_q(\Theta)=-i\omega\Theta^{\prime}(u)+\omega q \Theta(u)-\sum_{j=1}^{N_c} \kappa_j 
 {\tilde{\Delta}}_j\Theta(u).  
\end{equation}

We demonstrate that this mapping is invertible and get an upper  bound 
for the norm of its inverse. 

Applying the operator $M_q$ to the Fourier elements $\exp(i\,lu)$ in (\ref{eq:Fourier}) results in
\begin{equation}
 M_q\exp(i\,lu)=\nu_l(q) \exp(i\,lu),
\end{equation}
where 
\begin{equation}
 \nu_l(q)= \omega (q+l)+4\sum_{j=1}^{N_c}\,
 \kappa_j\sin^2\left(\frac{q+l}{2}j\right).
\end{equation}

By the assumption (\ref{eq:As}) one has $\nu_l(q) \ne 0$, $\forall \,\,l \in {\mathbb{Z}}$, 
so that the mapping $M_q$ possesses an inverse obeying 
$M_q^{-1}\exp(i\,lu)=(1/\nu_l)\exp(i\,lu)$. 
For the bounded 
linear operator $M_q^{-1}:\,\,X^{0}_q \rightarrow X_q^1$ one derives:
\begin{eqnarray}
 \parallel M_q^{-1} \parallel_{X^{0}_q,X^1_q}&\equiv& \parallel M_q^{-1}\parallel
=\sup_{0 \neq \Theta \in X^0_q}\frac{\parallel M^{-1}_q\,\Theta \parallel_{X^1_q}}
{\parallel \Theta \parallel_{X^0_q}}\nonumber\\
&=&\sup_{0 \neq \Theta \in X^0_q}\frac{\parallel \sum_{l\in {\mathbb{Z}}}\frac{1}{\nu_l}
\Theta_l \exp(i\,lu)\parallel_{X^1_q}}{\parallel \Theta \parallel_{X^0_q}}\nonumber\\
 &=& \sup_{0 \neq \Theta \in X^0_q} 
\frac{ \left( 
 \sup_{u\in [0,2\pi]}\left|\sum_{l\in {\mathbb{Z}}}\frac{1}{\nu_l} \Theta_l \exp(i\,lu)\right|
 +\sup_{u\in [0,2\pi]}
 \left|\left(\sum_{l\in {\mathbb{Z}}}\frac{1}{\nu_l}
 \Theta_l \exp(i\,lu)\right)^{\prime}\right|\right)}
 {\parallel \Theta \parallel_{X^0_q}}\nonumber\\
  &\le& \sup_{l \in {\mathbb{Z}}}\frac{1+|l|}{|\nu_l|} \,\sup_{0 \neq \Theta \in X^0_q} 
\frac{ 
 \sup_{u\in [0,2\pi]}\left|\sum_{l\in {\mathbb{Z}}} \Theta_l \exp(i\,lu)\right|}
 {\parallel \Theta \parallel_{X^0_q}}\nonumber\\
 &= & \sup_{l \in {\mathbb{Z}}} \frac{1+|l|}{|\nu_l|}\nonumber\\
 &=&\sup_{l\in {\mathbb{Z}}}\frac{1+|l|}{\left|(q+l)
 \left(\omega +\frac{4}{q+l}\sum_{j=1}^{N_c}\,
 \kappa_j\sin^2\left(\frac{q+l}{2}j\right)\right)\right|}\nonumber\\
&  \le& \left(\tilde{q}+\frac{1}{1-q}\right)\,\frac{1}{|\omega|-{\cal{R}}}\nonumber\\
& \le& \frac{(1+q){\cal R}+4{\bar{\kappa}}+
a(1+{\cal {P}}^b)}{a(1+{\cal {P}}^b){\cal R}},
\label{eq:boundL}
\end{eqnarray}
where $\tilde{q}$ is given in (\ref{eq:tildeq}).

For periodic travelling wave solutions $\Phi\in C^{1}_{2\pi}({\mathbb{R}},{\mathbb{C}})$) 
one derives, using (\ref{eq:A1}),(\ref{eq:norm}),(\ref{eq:As}) and (\ref{eq:ref}),  the bounds
\begin{equation}
\max_{u\in [0,2\pi]}|\Phi(u)|\le {\cal P}^{1/2},\,\,\,
\max_{u\in (0,2\pi)}|\Phi^{\prime}(u)|\le  \left(q+\frac{4{\bar{\kappa}}
 +a(1+{\cal P}^b)}{{\cal R}}\right){\cal P}^{1/2}.
\end{equation} 

We consider then the closed and convex subsets of $X_q^0$ and $X_q^1$ determined by 
\begin{equation}
 {{Y}}_q^0\,=\,\left\{\,\Theta\in X_q^0
 \,\,:
 \,\,\parallel \Theta\parallel_{C_{2\pi}^0}\le {\cal P}^{1/2}\,\right\}\,,
\end{equation}
and 
\begin{equation}
 {{Y}}_q^1\,=\,\left\{\,\Theta\in X_q^1
\,\,:
 \,\,\parallel \Theta\parallel_{C_{2\pi}^1}\le \left(1+q+\frac{4{\bar{\kappa}}
 +a(1+{\cal P}^b)}{{\cal R}}\right){\cal P}^{1/2}\,\right\}\,,
\end{equation}
respectively.  $Y_q^{1}$  is compactly embedded in  
$Y_q^0$ ($Y_q^1\Subset Y_q^{0}$).

Furthermore associated with the r.h.s. of (\ref{eq:ref1}) we introduce 
the nonlinear operator $N_q\,:\,Y_q^{0}\,\rightarrow\,Y_q^{0}$, as
\begin{equation}
 N_q(\Theta)= F(|\Theta|^2)\Theta.\label{eq:G}
\end{equation}
   
Clearly, the  operator $N_q$ 
is uniformly continuous on $Y_q^{0}$.  
The range is contained in a bounded ball in $Y_q^{0}$ 
since,
\begin{equation}
 \parallel N_q(\Theta) \parallel_{Y^{0}_q}=\parallel \, F(|\Theta|^2)\Theta\, \parallel_{C^0_{2\pi}}=
 \max_{u \in [0,2\pi]}|\, F(|\Theta(u)|^2)\Theta(u)|\,
 \le a(1+{\cal P}^{b})\,{\cal P}^{1/2}.
 \label{eq:rangeN}
\end{equation}

Finally, we express the problem  (\ref{eq:ref1})
 as a fixed point equation in terms of a mapping $Y_q^0\,\rightarrow\,Y_q^1\Subset Y_q^{0}$:
\begin{equation}
 \Phi=M_q^{-1}\,\circ\,N_q(\Phi)\equiv T_q(\Phi)\label{eq:compose}
\end{equation}

We get
\begin{equation}
 \parallel T_q(\Phi)\parallel_{Y_q^1}=\parallel M_q^{-1}(\,N_q(\Phi))\parallel
 \le \parallel M_q^{-1} \parallel \,\parallel  \,N_q(\Phi)\parallel_{Y_q^0}
 \le \left(1+q+\frac{4{\bar{\kappa}}
 +a(1+{\cal P}^b)}{{\cal R}}\right){\cal P}^{1/2},
\end{equation}
verifying that indeed
\begin{equation}
T_q(Y_q^0)\subseteq Y_q^1.
\end{equation}
Hence $T_q$ maps bounded subsets $Y_q^0$ of $X_q^0$ into relatively compact subsets 
$Y_q^1$ of  $Y_q^{0}$.

\vspace*{0.5cm}

It remains to prove that $T_q$ is continuous on $Y_q^{0}$. 
As $N$ is  uniformly continuous on $Y_q^{0}$,   
one has  $\forall t \in [0,2\pi]$ and $\forall \,\Phi_1,\Phi_2 \in Y_q^0$ that 
for a fixed arbitrary $\epsilon >0$ 
there exists $\delta >0$ such that 
\begin{equation*}
 \parallel N_q(\Phi_1)(t)-N_q(\Phi_2)(t)\parallel_{Y_q^{0}}<
\frac{a(1+{\cal {P}}^b){\cal R}}{(1+q){\cal R}+4{\bar{\kappa}}+
a(1+{\cal {P}}^b)}\,\epsilon
\end{equation*}
if $\parallel \Phi_1-\Phi_2\parallel_{Y_q^0} < \delta$.
Hence, for arbitrary $\Phi_1,\Phi_2\in Y_q^{0}$ we have
\begin{equation*}
 \parallel M^{-1}_q(N_q(\Phi_1))- M^{-1}_q(N_q(\Phi_2))
 \parallel_{Y_q^1\Subset Y_q^{0}}\le 
 \parallel M^{-1}_q\parallel
 \,\parallel N_q(\Phi_1)- N_q(\Phi_2)\parallel_{Y_q^{0}}<\epsilon,
\end{equation*}
verifying g that $T_q(\Phi)$ is continuous on $Y_q^0$.
 Schauder’s fixed point theorem
implies then that the fixed point equation $\Phi =  T_q(\Phi)$ has at least one solution.

Furthermore, the spectrum of linear plane wave solutions (phonons) arising for zero nonlinear term,
determined by the r.h.s. of  
system (\ref{eq:ref1}), 
forms a continuous band with values in the interval $[-{\cal R},{\cal R}]$. 
However,  since by the hypothesis  (\ref{eq:As}) the values of the 
frequency of oscillations $\omega$ lie outside the range of the linear (phonon) band 
the corresponding orbits 
are anharmonic. This necessitates
amplitude-depending  tuning of the frequency so that the latter comes to lie 
outside of the phonon spectrum. 
Thus it must hold that 
$\parallel N_q(\Phi)\parallel_{Y^0_q} =\parallel \, F(|\Phi|^2)\Phi\, \parallel_{C^0_{2\pi}} \not \equiv 0$  
which is fulfilled only if  $\Phi \not \equiv 0$. That is, the fixed point equation (\ref{eq:compose}) possesses only nonzero 
solutions   and the proof is finished.

\vspace*{0.5cm}

\hspace{16.5cm} $\square$

\vspace*{1.0cm}

To summarise, we have proven the existence  of nonzero periodic 
travelling wave solutions 
for a general DNLS (including as a special case the standard DNLS) on finite one-dimensional lattices. To this end the existence problem has been 
reformulated  
as a fixed point problem for an operator on a function space  
which is solved with the help  of Schauder's Fixed Point Theorem. 
Our method can be straightforwardly extended to treat also the general DNLS on lattices of 
higher dimension.

\end{document}